\documentclass[a4paper,11pt]{article}
\usepackage{amssymb}
\usepackage{amsmath}
\usepackage{amsthm}
\usepackage{graphicx}
\usepackage{cite}
\usepackage[square,numbers,sort&compress]{natbib}
\usepackage[utf8]{inputenc}
\newtheorem{thm}{Theorem}[section]
\newtheorem{cor}{Corollary}[section]
\newtheorem{lem}{Lemma}[section]

\theoremstyle{definition}

\theoremstyle{remark}
\newtheorem{rem}{Remark}[section]

\numberwithin{equation}{section}

\begin{document}
	
	\begin{center}
		{\Large \bf Bounds on Derivatives in Compositions of Two Rational Functions with Prescribed Poles }
	\end{center}
	\begin{center}
		{\normalsize \bf Preeti Gupta$^{(a)}$}
	\end{center}
	\begin{center}
		
		{\normalsize $^{(a)}$Department of Applied Mathematics\\ Amity University, Noida-201313, India.\\ Email: preity8315@gmail.com}

	\end{center}
	
	\begin{abstract}
This paper explores a class of rational functions \( r(s(z)) \) with degree \( mn \), where \( s(z) \) is a polynomial of degree \( m \). Inequalities are derived for rational functions with specified poles, extending and refining previous results in the field.
	\end{abstract}

	\noindent\textbf{Mathematics Subject Classification (2010)}: 30A10, 30C15, 26D07.\\
	\noindent\textbf{Keywords and phrases}: Rational functions; Polynomial inequalities; Zeros; Poles; Blaschke product.

	\section{INTRODUCTION}	

Let \( P_n \) represent the class of all complex polynomials of degree no more than \( n \), and let \( k \) be a positive real number. Define the following regions in the complex plane:
\[
T_k = \{ z : |z| = k \}, \quad D_k^- = \{ z : |z| < k \}, \quad D_k^+ = \{ z : |z| > k \}.
\]
Consider a polynomial \( p(z) \) of degree \( n \). In 1926, Bernstein \cite{Bernstein1912} proved the following well-known inequality:
\[
\max_{|z| = 1} |p'(z)| \leq n \max_{|z| = 1} |p(z)|. \tag{1}
\]
Equality in \((1)\) holds only when \( p(z) = az^n \), where \( a \neq 0 \).\\

If we restrict to polynomials with no zeros in \( D_1^- \), inequality \((1)\) can be improved. P. Erdős conjectured and Lax \cite{Lax1944} later proved that for such polynomials \( p(z) \):
\[
\max_{|z| = 1} |p'(z)| \leq \frac{n}{2} \max_{|z| = 1} |p(z)|. \tag{2}
\]

For polynomials with no zeros in \( D_1^+ \), Turán \cite{PT} established the following inequality:
\[
\max_{|z| = 1} |p'(z)| \geq \frac{n}{2} \max_{|z| = 1} |p(z)|. \tag{3}
\]\\

For $a_{j}\in\mathbb{C} $ with $j=1,2,\dots,n,$ let\\

$w(z):=\prod_{j=1}^{n}\left(z-a_{j} \right), $\hspace{.3cm}	$B(z):=\frac{w^{*}(z)}{w(z)}=\prod_{j=1}^{n}\left(\frac{1-\bar{a}_{j}z}{z-a_{j}} \right),$\\
and
$$ \Re_{n}=\Re_{n}\left(a_{1},a_{2},\dots,a_{n} \right):=\left\lbrace\frac{p(z)}{w(z)}:p\in P_{n} \right\rbrace  . $$
The product $B(z)$ is known as Blaschke product and one can easily verify that $ \left|B(z) \right| =1 $ and $\frac{zB^{'}(z)}{B(z)}=|B^{'}(z)| $ for $ z\in T_{1} $. Then $\Re_{n}$ is the set of all rational functions with at most $n$ poles $ a_{1},a_{2},\dots,a_{n} $ and with finite limit at infinity. We observe that $ B(z)\in \Re_{n} $. For $f$ defined on $T_{1}$ in the complex plane, we denote 
$\left \| f \right \|=\sup_{z\in T_{1}} {|f(z)|},$ the Chebyshev norm of $f$ on $T_{1}.$
Throughout this paper, we always assume that all poles $ a_{1},a_{2},\dots,a_{n} $ are in $D_{1}^{+}.$\\
In 1995, Li, Mohapatra and Rodriguez~\cite{Li1995} have proved Bernstein-type  inequalities for rational function $ r(z)\in \Re_{n} $ with prescribed poles where they replaced $z^{n}$ by Blaschke product $B\left(z \right) $ and established the following results.\\ 

\begin{thm}{\label{thm1.1}}	
	If $ r\in \Re_{n} $, and all zeros of r lie in $T_{1}\cup D_{1}^{+},$ then for $z \in T_{1}$, we have 
	\begin{align}{\label{1.1}}
		\left|{r}'(z) \right| \leq\frac{1}{2}\left|{B}'(z) \right|\left\| r\right\|. \tag{4}
	\end{align}
	Equality holds for $r(z)=aB(z)+b$ with $\left|a \right| =\left|b \right| =1.$
\end{thm}
\begin{thm}{\label{thm1.2}}
	Let $ r\in \Re_{n}, $ where r has exactly n poles at $ a_{1},a_{2},\dots,a_{n} $ and all its zeros lie in $T_{1}\cup D_{1}^{-}.$ Then for $z\in T_{1}$
	\begin{align}{\label{1.2}}
		\left|{r}'(z) \right| \geq \frac{1}{2}\left[\left|{B}'(z) \right|-\left( n-t\right)   \right] \left|r(z) \right|, \tag{5}
	\end{align}  
	where t is the number of zeros of r with counting multiplicity. The above result is best possible and equality holds for $r(z)=aB(z)+b$ with $\left|a \right| =\left|b \right| =1.$ 	
\end{thm}
\begin{rem}
	In particular, if $r$ has exactly $n$ zeros in $T_{1}\cup D_{1}^{-},$ then the inequality (\ref{1.2}) yields Bernstein-type inequality, namely for $z \in T_{1}$
	\begin{align}{\label{1.3}}
		\left|{r}'(z) \right| \geq \frac{1}{2}\left|{B}'(z) \right|\left|r(z) \right|. \tag{6}	
	\end{align}
\end{rem}
Aziz and Shah~\cite{Aziz1997} proved the following theorems which improves upon the inequalities (\ref{1.2}) and (\ref{1.3}) by introducing $m=\min_{z\in T_{1}}\left|r(z) \right|.$
\begin{thm}{\label{thm1.3}}
	If $ r\in \Re_{n} $, with all its zeros lie in $T_{1}\cup D_{1}^{+},$ then for $z \in T_{1},$ we have 
	\begin{align}{\label{1.4}}
		\left|{r}'(z) \right| \leq\frac{1}{2}\left|{B}'(z) \right|\left( {\left \| r \right \|}-m \right), \tag{7}
	\end{align} 
	where $m=\min_{z\in T_{1}}\left|r(z) \right|.$ The result is best possible and equality attains for $r(z)=B(z)+he^{i\alpha }$ with $h\geq1$ and $\alpha$ real.
\end{thm}
\begin{thm}{\label{thm1.4}}
	Let $ r\in \Re_{n}, $ where r has exactly n poles at $ a_{1},a_{2},\dots,a_{n} $ and all its zeros lie in $T_{1}\cup D_{1}^{-}.$ Then for $z\in T_{1}$
	\begin{align}{\label{1.5}}
		\left|{r}'(z) \right| \geq\frac{1}{2}\left|{B}'(z) \right|\left[ \left|r\left(z \right) \right| +m \right] , \tag{8}
	\end{align} 
	where $m=\min_{z\in T_{1}}\left|r(z) \right|.$ 
	Equality attains for $r(z)=B(z)+he^{i\alpha }$ with $h\leq1$ and $\alpha$ real.	
\end{thm}

Arunrat and Nakprasit \cite{Arunrat2020} proved the following results, which not only improve upon the inequalities (\ref{1.4}) and (\ref{1.5}), but also generalize them.

\begin{thm}{\label{thm1.5}}
	Let $r \in \Re_{n},$ where r has exactly n poles at $ a_{1},a_{2},\dots,a_{n} $ and all its zeros lie in $T_{k}\cup D_{k}^{+},$ $k\geq 1.$ Then for $z\in T_{1}$
	\begin{align}
		\left|{r}'(z) \right|\leq \frac{1}{2}\left[\left|{B}'(z) \right|-\frac{\left( n(1+k)-2t\right)\left(\left|r(z) \right|  -m\right)^{2}  }{\left( 1+k\right)\left( \left \| r  \right \| -m \right)^{2}  }  \right]\left(\left \| r  \right \| -m \right),   \tag{9}
	\end{align}
	where t is the number of zeros of r with counting multiplicity and $m=\min_{z\in T_{k}}\left|r(z) \right|$.
\end{thm}

\begin{thm}{\label{thm1.6}}
	Let $ r\in \Re_{n} $, where r has exactly n poles $ a_{1},a_{2},\dots,a_{n} $ and  all its zeros lie in $T_{k}\cup D_{k}^{-}$, $k \leq 1$. Then for $z \in T_{1}$ 
	\begin{align}
		\left|{r}'(z) \right| \geq \frac{1}{2}\left[\left|{B}'(z) \right|+\frac{2t-n(1+k)}{1+k} \right]\left( \left|r(z) \right|+m\right) , \tag{10}	
	\end{align}
	where $t$ is the number of zeros of $r$ with counting multiplicity and $m=\min_{z\in T_{k}}\left|r(z) \right|$.
\end{thm}

In 2015, Qasim and Liman \cite{Qasim} investigated a class of rational functions, \( r(s(z)) \in \Re_{mn} \), where all the poles \( a_1, a_2, \dots, a_{mn} \) are contained within the region \( D_1^+ \). These rational functions are expressed as the composition \( r(s(z)) \), defined by  
\[
(r \circ s)(z) = r(s(z)) =: \frac{p(s(z))}{w(s(z))}, \tag{11}
\]  
where \( s(z) \) is a polynomial of degree \( m \) with all its zeros confined to \( T_1 \cup D_1^- \), and \( r \) belongs to \( \Re_n \).  

The term \( w(s(z)) \) represents the product of factors corresponding to the poles and is given by  
\[
w(s(z)) = \prod_{j=1}^{mn} (z - a_j). \tag{12}
\]  

Additionally, the Blaschke product associated with \( w(s(z)) \) is defined as  
\[
B(z) = \frac{w^*(s(z))}{w(s(z))} = \frac{z^{mn} w(s(1/\bar{z}))}{w(s(z))} = \prod_{j=1}^{mn} \frac{1 - \bar{a}_j z}{z - a_j}. \tag{13}
\]  
Qasim and Liman \cite{Qasim} provided the following generalization of a known inequality.  

\begin{thm}{\label{thm1.7}}
Suppose \( r(s(z)) \in \Re_{mn} \), where \( r(s(z)) \) has no zeros in \( D_1^- \), and all zeros of \( s(z) \) lie within \( T_1 \cup D_1^- \). Then for all \( z \in T_1 \), we have  
\[
\left| r'(s(z)) \right| \leq \frac{1}{2mm'} \left| B'(z) \right| \cdot \| r \circ s \|, \tag{14}
\]  
where \( m' = \min_{z \in T_1} |s(z)| \). The inequality is sharp, and equality holds for \( r(s(z)) = aB(z) + b \), where \( a, b \in T_1 \) and \( s(z) = z^m \). 	
\end{thm}

Note that if \( s(z) \) has a zero on \( T_1 \), then \( m' = 0 \), resulting in a trivial inequality:  
\[
0 \leq 2mm' \cdot \left| r'(s(z)) \right| \leq \left| B'(z) \right| \cdot \| r \circ s \|. \tag{15}
\]

\begin{thm}{\label{thm1.8}}
Let \( r(s(z)) \in \Re_{m,n} \) with \( r(s(z)) \neq 0 \) in \( D_1^{-} \), and let all zeros of \( s(z) \) lie in \( T_1 \cup D_1^{-} \). Then for \( z \in T_1 \), we have

\[
\left| r'(s(z)) \right| \leq \frac{1}{2mm'} \left| B'(z) \right| \cdot\left(  \| r \circ s \| - m^*\right) , \tag{16}
\] 

where
\[
m' = \min_{z \in T_1} |s(z)|, \quad m^* = \min_{z \in T_1} |r(s(z))|.
\]
Equality holds for \( r(s(z)) = B(z) + h e^{i\alpha} \), where \( s(z) = z^m \), \( h \geq 1 \), and \( \alpha \) is real.
\end{thm}

Arunrat and Nakprasit \cite{ArunratN2020} generalized the Theroem \ref{thm1.8}.

\begin{thm}{\label{thm1.9}}
	
	Let \( r(s(z)) \in \Re_{m,n} \) with \( r(s(z)) \neq 0 \) in \( D_k^{-} \), \( k \geq 1 \), and let all zeros of \( s(z) \) lie in \( T_1 \cup D_1^{-} \). Then for \( z \in T_1 \), we have
	
\begin{equation}
	\left| r'(s(z)) \right| \leq \frac{1}{2mm'}\left[  \left| B'(z) \right| - \frac{(mn(1 + k) - 2mt) (|r(s(z))| - m^*)^{2}}{(1 + k) (\|r \circ s \| - m^*)^{2}}\right] 
	 \left(  \|r \circ s\| - m^*\right) , \tag{17}	
\end{equation}	
where $	mt$ is the number of zeros of $r \circ s $ with counting multiplicity,
	
	\[
	m' = \min_{z \in T_1} |s(z)|, \quad m^* = \min_{z \in T_k} |r(s(z))|.
	\]
Equality holds for 
	\[
	r(s(z)) = \frac{(z + k)^{mt}}{(z - a)^{mn}}, \quad s(z) = z^m \text{and} \quad B(z) = \left(\frac{1 - az}{z - a}\right)^{mn},
	\]
	where \( a > 1 \), \( k \geq 1 \), at \( z = 1 \).
	
\end{thm}


	\section{MAIN RESULTS}	
	
In this paper, we consider the class of rational functions \( \Re_{mn} \) having no zeros in \( D_k^- \), where \( k \geq 1 \), and establish a generalization of the result by Arunrat and Nakprasit \cite{ArunratN2020}.
We shall obtain bounds for the derivative of composition of two rational functions by involving the moduli of all its zeros. More precisely, we have the following.

	\begin{thm}{\label{thm2.1}}
		Let \( r(s(z))= \frac{p\left( s(z)\right) }{w\left( s(z)\right) } \in \Re_{mn} \) with \( r(s(z)) \neq 0 \) in \( D_k^{-} \), \( k \geq 1 \), and let all zeros of \( s(z) \) lie in \( T_1 \cup D_1^{-} \). Then for \( z \in T_1 \),
	we have
	\begin{equation*}{\label{2.1}}
\left|{r }'\left( s(z) \right) \right|\leq \frac{1}{2mm^{'}}\left[\left|{B}'(z) \right|-\frac{2\left[  \frac{mn}{2}- \left( \sum_{j=1}^{mt} \frac{1}{1+\left|b_{j} \right| }\right) \right]   \left(\left|r(s(z)) \right| -m^{*} \right)^{2}}{\left(\left \| r \circ s \right \|-m^{*} \right)^{2}}  \right] \left(\left \| r \circ s \right \|-m^{*} \right),	\tag{18}	
	\end{equation*}
where $mt$ is the number of zeros of $ r \circ s $ with counting multiplicity,
	\[
	m' = \min_{z \in T_1} |s(z)|, \quad m^* = \min_{z \in T_k} |r(s(z))|.
	\]
Equality holds for 
	\[
	r(s(z)) = \frac{(z + k)^{mt}}{(z - a)^{mn}}, \quad s(z) = z^m, \quad B(z) = \left(\frac{1 - az}{z - a}\right)^{mn},
	\]
	where \( a > 1 \), \( k \geq 1 \), at \( z = 1 \).
	
\end{thm}

\begin{rem}As we have
		\begin{equation*}{\label{re3.4}}
			\frac{1}{1+\left|b_{j} \right| } \leq \frac{1}{k+1}, \tag{19}
		\end{equation*}
		for $\left|b_{j} \right|\geq k \geq 1 .$ Using (\ref{re3.4}) in (\ref{thm2.1}), we observe that   inequality ({\ref{thm2.1}}) reduces to the Theorem \ref{thm1.9}.	
	\end{rem}
	
	\begin{thm}{\label{thm2.2}}
		Let $r\left( s(z) \right)= \frac{p\left( s(z)\right) }{w\left( s(z)\right) } \in \Re_{mn}.$ If $b_{1},b_{2},\dots,b_{mt}$ are the zeros of $r\left( s(z) \right)$ all lying in $T_{k}\cup D_{k}^{-},$ $k\leq 1$. Then for $z \in T_{1},$ we have 
		\begin{equation}{\label{2.3}}
	 \left|{r^{'}(s(z))}  \right| \geq  \frac{1}{2mm^{'}}\left[ \left|{B}'\left(z \right)  \right| + 2\left(\sum_{j=1}^{mt}\frac{1}{1+\left|b_{j} \right| }-\frac{mn}{2} \right)  \right]\left( \left|r\left(s(z) \right)\right|  + m^{*} \right).  \tag{20}			
		\end{equation}
Where $m^{*}=\min_{z\in T_{k}}\left|r(s(z)) \right|$.
	\end{thm}
	\begin{rem}As
		\begin{equation}{\label{re3.2}}
			\frac{1}{1+\left|b_{j} \right| } \geq \frac{1}{1+k}, \tag{21}
		\end{equation}
		where $\left|b_{j} \right|\leq k \leq 1 .$ Using (\ref{re3.2}) in (\ref{2.3}), we see that inequality reduces to the following result.
	\end{rem}
		\begin{align*}
		\left|{r}'(s(z)) \right| \geq \frac{1}{2mm^{'}}\left[\left|{B}'(z) \right|+\frac{2mt-mn(1+k)}{1+k} \right]\left( \left|r(s(z)) \right|+m^{*}\right) . 	
	\end{align*}
Specifically, setting \( k = 1 \) and \( t = n \), we obtain the subsequent result.

\begin{cor}
Consider \( r\left( s(z) \right) = \frac{p\left( s(z)\right)}{w\left( s(z)\right)} \in \Re_{mn}. \) Suppose \( b_{1}, b_{2}, \dots, b_{mn} \) are the zeros of \( r\left( s(z) \right) \) all contained in \( T_{1} \cup D_{1}^{-} \). Then, we have 
\begin{equation*}
	\left|{r^{'}(s(z))}  \right| \geq  \frac{1}{2mm^{'}} \left|{B}'\left(z \right)\right| \left[  \left|r\left(s(z) \right)\right|  + m^{*} \right] . 			
\end{equation*}
Where $m^{*}=\min_{z\in T_{k}}\left|r(s(z)) \right|$.
\end{cor} 	
By taking $s(z)=z$ in Theorems \ref{thm2.1} and \ref{thm2.2}, the results simplify to those obtained by Preeti Gupta et al. \cite{GHM}. This is because substituting $s(z)=z$ reduces the composite function $r(s(z))$ to $r(z)$, thereby recovering the specific results established in their work as a special case of these more general theorems.
	\section{LEMMAS}		
	
For the proof of these theorems, we need the following lemmas. The first lemma is due to Li, Mohapatra and Rodriguez~\cite{Li1995}.

\begin{lem}{\label{lem 2.1}}
	If $r \in \Re_{n}$ and $r^{*}\left ( z \right )= B\left ( z \right )\overline{r\left ( \frac{1}{\bar{z}} \right )},$ then for $z\in T_{1}$ 
	\begin{align*}
		\left|{\left(r^{*}(z) \right)}'  \right|+\left|{r}'(z) \right| \leq \left|{B}'(z) \right|\left \| r \right \|.   
	\end{align*}
\end{lem}
Lemma 2 is due to Aziz and Dawood \cite{AQ}.
\begin{lem}{\label{lem 2.2}}
	If \( p \in P_n \) and \( p(z) \) has its zeros in \( T_1 \cup D_1^- \), then
	\[
	\min_{z \in T_1} |p'(z)| \geq n \cdot \min_{z \in T_1} |p(z)|. 
	\]
	The inequality is sharp, and equality holds for polynomials having all zeros at the origin.
\end{lem}

		\section{PROOF OF THEOREMS}	
	
	\begin{proof}[\textbf{Proof of Theorem 2.1}]Assume that $r(s(z)) \in \Re_{mn}$ has no zero in $D^{-}_{k}$, where $k \geq 1$. Let 
		\[
		m^{*} = \min_{z \in T_{k}} \left| r(s(z)) \right|
		\]
		and let \( tm \) be the number of zeros of \( r(s(z)) \), counting multiplicities. Therefore, 
		\[
		m^{*} \leq \left| r(s(z)) \right| \quad \text{for } z \in T_{k}.
		\]		
If \( r(s(z)) \) has a zero on \( T_{k} \), then \( m^{*} = 0 \), and hence for every \( \alpha \) with \( \left| \alpha \right| < 1 \), we get 
		\[
		r(s(z)) - \alpha m^{*} = r(s(z)).
		\]
		In case \( r(s(z)) \) has no zero on \( \left| z \right| = k \), we have for every \( \alpha \) with \( \left| \alpha \right| < 1 \) that 
		\[
		\left| -\alpha m^{*} \right| = \left| \alpha \right| m^{*} \leq \left| r(s(z)) \right| \quad \text{for } \left| z \right| = k.
		\]
By Rouch\'e's theorem, 
		\[
		R(z) = r(s(z)) - \alpha m^{*}
		\]
		and \( r(s(z)) \) have the same number of zeros in \( \left| z \right| < k \). That is, for every \( \alpha \) with \( \left| \alpha \right| < 1 \), \( R(z) \) has no zeros in \( \left| z \right| < k \).
		
		If \( p(z) \) has \( t \) zeros and \( s(z) \) has \( m \) zeros, then \( p(s(z)) \) has \( tm \) zeros. Let \( b_1, b_2, \dots, b_{mt} \) be the zeros of \( p(s(z)) \), where \( mt \leq mn \). Now,
		\[
		r(s(z)) = \frac{p(s(z))}{w(s(z))}.
		\]
		\begin{align*}
			\frac{z{R}'\left(z \right) }{R\left(z \right) }
			=  \frac{z({r(s(z))})'}{r(s(z)) }
			&=	 \frac{z({p(s(z))})' }{p(s(z)) } - \frac{z{w(s(z))}'}{w\left(s(z) \right) }\\
			&= \sum_{j=1}^{mt} \frac{z}{z - b_j} - \sum_{j=1}^{mn} \frac{z}{z - a_j}.
		\end{align*}
	Then
	\begin{align*}
		Re\left(\frac{z{R}'\left(z \right) }{R\left(z \right) } \right) 
	= Re \left(\sum_{j=1}^{mt}\frac{z}{z-b_{j}} \right)  - Re\left(\sum_{j=1}^{mn} \frac{z}{z - a_j} \right) .
	\end{align*}
		\begin{align*}{\label{4.1}}
			Re\left(\frac{z{R}'\left(z \right) }{R\left(z \right) } \right) 
			&\leq \sum_{j=1}^{mt}\frac{1}{1+\left|b_{j} \right| }  + \sum_{j=1}^{mn} Re\left( \frac{1}{2}-\frac{z}{z - a_j}\right) -\frac{mn}{2}  \nonumber\\
			& = \sum_{j=1}^{mt}\frac{1}{1+\left|b_{j} \right| } -\frac{mn}{2}+ \sum_{j=1}^{mn} \frac{|a_j|^2 - 1}{2|z - a_j|^2}  \nonumber\\
			& = \sum_{j=1}^{mt}\frac{1}{1+\left|b_{j} \right| } -\frac{mn}{2}+\frac{|B^{'}(z)|}{2}. \tag{22} \nonumber\\
		\end{align*}
		Note that $R^{*}\left(z \right)=B\left(z \right)\overline{R\left( \frac{1}{\bar{z}}\right) } = B\left(z \right) \bar{R}\left(\frac{1}{z} \right).$ Then 
		\begin{align*}
			{\left( R^{*}\left(z \right)\right) }'
			&= {B}'\left(z \right)\bar{R}\left(\frac{1}{z} \right) + B\left(z \right) \left( {\bar{R}\left(\frac{1}{z} \right)}'\right)  \\
			&={B}'\left(z \right) \bar{R}\left(\frac{1}{z} \right) +  B\left(z \right) \left( {\bar{R}}'\left(\frac{1}{z} \right)\right)  \left(-\frac{1}{z^{2}} \right) \\
			&={B}'\left(z \right) \bar{R}\left(\frac{1}{z} \right) - \frac{B\left( z\right) }{z^{2}}\left( {\bar{R}}'\left(\frac{1}{z} \right)\right). 
		\end{align*}
	This implies 	 
\begin{align*}
	 z{\left( R^{*}\left(z \right)\right) }'= z{B}'\left(z \right) \bar{R}\left(\frac{1}{z} \right) - \frac{B\left( z\right) }{z}\left( {\bar{R}}'\left(\frac{1}{z} \right)\right).
\end{align*}	
		
		Since $z \in T_{1},$ we have $\bar{z}=\frac{1}{z},$
		$\left|B\left(z \right)  \right|=1, \frac{z{B}'\left(z \right) } {B\left( z\right) }=\left|{B}'\left(z \right)  \right|, $ hence 
		\begin{align*}
			\left|z{\left( R^{*}\left(z \right)\right) }' \right|
			&= \left|z{B}'\left(z \right)\overline{R\left(z \right) }  - B\left( z\right)\overline{z{R}'\left(z \right) }  \right| \\
			&= \left|\frac{z{B}'\left(z \right)}{B\left(z \right) }.\overline{R\left(z \right) }  - \overline{z{R}'\left(z \right) }  \right| \\
			&=\left|\left|{B}'\left(z \right)  \right|\overline{R\left(z \right) } -\overline{z {R}'\left(z \right)}  \right|. 	
		\end{align*}
		Since $\left|{B}'\left(z \right)  \right|$ is real, we get\\
		
		$\left|z{\left( R^{*}\left(z \right)\right) }' \right|=\left|\left|{B}'\left(z \right)  \right|R\left(z \right)  -z {R}'\left(z \right)  \right|.$
		Then by inequality (\ref{4.1})
		
		\begin{align*}
			\left|\frac{z{\left( R^{*}\left(z \right)\right) }'}{R\left(z \right) } \right|^{2}
			&=\left|\left|{B}'\left(z \right)  \right|-\frac{z{R}'\left(z \right) }{R\left(z \right) } \right| ^{2}\\
			&=\left|{B}'\left(z \right)  \right|^{2}+\left| \frac{z{R}'\left(z \right) }{R\left(z \right) }\right| ^{2}-2\left|{B}'\left(z \right)  \right|Re\left(\frac{z{R}'\left(z \right) }{R\left(z \right) } \right) \\
			& \geq \left|{B}'\left(z \right)  \right|^{2}+\left| \frac{z{R}'\left(z \right) }{R\left(z \right) }\right| ^{2}\\
			&-2\left|{B}'\left(z \right)  \right|\left[\frac{|B^{'}(z)|}{2}+ \left(  \sum_{j=1}^{mt}\frac{1}{1+\left|b_{j} \right| } -\frac{mn}{2}\right)  \right] \\
			& =\left| \frac{z{R}'\left(z \right) }{R\left(z \right) }\right| ^{2}+ 2 \left[\frac{mn}{2}- \left(\sum_{j=1}^{mt}\frac{1}{1+\left|b_{j} \right| } \right) \right]\left|{B}'\left(z \right)  \right|. 
		\end{align*}
	This implies that for $z \in T_{1}$
		\begin{align}{\label{4.2}}
			\left[\left|{R}'\left(z \right)  \right|^{2} + 2\left[\frac{mn}{2}- \left(\sum_{j=1}^{mt}\frac{1}{1+\left|b_{j} \right| } \right)\right] \left|{B}'\left(z \right)  \right|\left|R\left(z \right)  \right| ^{2} \right] ^{\frac{1}{2}}\leq \left|{\left( R^{*}\left(z \right)\right) }' \right| ,	\tag{23}
		\end{align}
		where 
		$R^{*}\left(z \right)=B\left(z \right)\overline{R\left(\frac{1}{\bar{z}} \right) }= r^{*}\left(s(z) \right)-\bar{\alpha}m^{*} B\left(z \right) .$\\
		
		Moreover, $$ {\left(  R^{*}\left(z \right)\right) }' =\ {\left( r^{*}\left(s(z) \right)\right) }' -\bar{\alpha}m^{*} {B}'\left(z \right)$$
		and
		$${R}'\left(z \right)=\left( {r}\left(s(z) \right)\right) ^{'} ={\left(r\left( s(z)\right)-\alpha m^{*}  \right) }'.$$
		
Apply these relations into (\ref{4.2}), we obtain that  
$$
\left[\left|{(r(s(z)))}'  \right|^{2} +2\left[ \frac{mn}{2}-\left( \sum_{j=1}^{mt}\frac{1}{1+\left|b_{j} \right| } \right)\right] \left|{B}'\left(z \right)  \right|\left|r\left(s(z) \right)-\alpha m^{*}  \right| ^{2} \right] ^{\frac{1}{2}}\hspace{2.5cm}$$
\begin{equation}{\label{4.3}}
\hspace{7.3cm}\leq \left|{\left( r^{*}\left(s(z) \right)\right) }'- \bar{\alpha}m^{*}{B}'\left(z \right)  \right|,	\tag{24}
\end{equation}
 for $z \in T_{1}$ and for $\alpha$ with $\left|\alpha \right|<1. $
		Choose the argument of $\alpha$ such that 
		$\left|{(r^{*}\left(s(z) \right)) }'-\bar{\alpha}m^{*} {B}'\left(z \right)  \right|=\left|{\left( r^{*}\left(s(z) \right)\right)  }' \right|-m^{*}\left|\alpha \right|\left|{B}'\left(z \right)  \right|,   $
		for $z \in T_{1}.$\\
			Since  $$ \left|r\left(s(z) \right)-m^{*}\alpha  \right| \geq \left|\left|r\left(s(z) \right)  \right|-m^{*}\left|\alpha \right|   \right|. $$
		 Note that  $$ \left|  \left|r\left(s(z) \right)\right| -m^{*} \left| \alpha\right| \right| ^{2} =\left( \left|r\left(s(z) \right)  \right|-m^{*}\left|\alpha \right|\right)^{2},$$ 
		which implies that
		$$ \left|r\left(s(z) \right)-m^{*}\alpha  \right|^{2} \geq \left( \left|r\left(s(z) \right)  \right|-m^{*}\left|\alpha \right|  \right) ^{2}.$$
		Apply above inequality in (\ref{4.3})
$$	\left[\left|{\left( r(s(z))\right) }'  \right|^{2} +2\left[\frac{mn}{2}- \left( \sum_{j=1}^{mt}\frac{1}{1+\left|b_{j} \right| } \right)\right] \left|{B}'\left(z \right)  \right|\left( \left| r\left(s(z) \right)\right| - m^{*}\left| \alpha \right|  \right)  ^{2} \right] ^{\frac{1}{2}}\hspace{2.5cm}$$
\begin{equation*}
\hspace{7.3cm}\leq\left|{\left( r^{*}\left(s(z) \right)\right)  }' \right|-m^{*}\left|\alpha \right|\left|{B}'\left(z \right)  \right|.	
\end{equation*}
Letting $\left|\alpha \right|\rightarrow 1 $, we get 
$$	\left[\left|{\left( r(s(z))\right) }'  \right|^{2} +2\left[\frac{mn}{2}- \left( \sum_{j=1}^{mt}\frac{1}{1+\left|b_{j} \right| } \right)\right] \left|{B}'\left(z \right)  \right|\left( \left| r\left(s(z) \right)\right| - m^{*}  \right)  ^{2} \right] ^{\frac{1}{2}}\hspace{2.3cm}$$
\begin{equation*}
\hspace{7.3cm}\leq\left|{\left( r^{*}\left(s(z) \right)\right)  }' \right|-m^{*}\left|{B}'\left(z \right)  \right|.
\end{equation*}
By lemma \ref{lem 2.1}, implies that 
$$\left[\left|{\left( r(s(z))\right) }'  \right|^{2} +2\left[\frac{mn}{2}- \left( \sum_{j=1}^{mt}\frac{1}{1+\left|b_{j} \right| } \right)\right] \left|{B}'\left(z \right)  \right|\left( \left| r\left(s(z) \right)\right| - m^{*}  \right)  ^{2} \right] ^{\frac{1}{2}}\hspace{2.3cm}$$
\begin{equation*}
\hspace{6.7cm}\leq \left|{B}'\left(z \right)  \right|\left \| r\circ s \right \|-\left|{\left( r(s(z))\right) }'  \right| -m^{*} \left|{B}'\left(z \right)  \right|. 	
\end{equation*}
Further simplifying and square both sides, gives us	
$$\left|{\left( r(s(z))\right) }' \right|^{2} +2\left[\frac{mn}{2}- \left( \sum_{j=1}^{mt}\frac{1}{1+\left|b_{j} \right| } \right)\right] \left|{B}'\left(z \right)  \right|\left( \left| r\left(s(z) \right)\right| - m^{*}  \right)  ^{2}\hspace{2.3cm}$$
\begin{equation*}
\leq  \left(\left \| r \circ s \right \|-m^{*}\right)^{2}\left|{B}'\left(z \right)  \right|^{2} +\left|{\left( r(s(z))\right) }' \right|^{2}\\
- 2\left(\left \| r \circ s \right \|-m^{*}\right)\left|{B}'\left(z \right)  \right|\left|{\left( r(s(z))\right) }'  \right|.
\end{equation*}
	This implies that
		\begin{align*}
			\left|{\left( r(s(z))\right) }'  \right| \leq \frac{\left(\left \| r \circ s \right \|-m^{*}\right)^{2}\left|{B}'\left(z \right)  \right|^{2}}{ 2\left(\left \| r \circ s \right \|-m^{*}\right)\left|{B}'\left(z \right)  \right|}
			-\frac{2\left[\frac{mn}{2}-\left( \sum_{j=1}^{mt}\frac{1}{1+\left|b_{j} \right| }\right) \right]  \left|{B}'\left(z \right)  \right|\left(\left|r\left(s(z) \right)  \right| -m^{*} \right)^{2}  }{2\left(\left \| r \circ s \right \|-m^{*}\right)\left|{B}'\left(z \right)  \right|}.
		\end{align*}
Thus 
\begin{equation}{\label{4.4}}
		\left|{\left( r(s(z))\right) }' \right|\leq \frac{1}{2}\left[\left|{B}'(z) \right|-\frac{2\left[  \frac{mn}{2}- \left( \sum_{j=1}^{mt} \frac{1}{1+\left|b_{j} \right| }\right) \right]   \left(\left|r(s(z)) \right| -m^{*} \right)^{2}}{\left(\left \| r \circ s \right \|-m^{*} \right)^{2}}  \right] \left(\left \| r \circ s \right \|-m^{*} \right). \tag{25}
\end{equation}
 For \( z \in T_1 \) with \( R(z) \neq 0 \), we have
 \[
 (r(s(z)))' = r'(s(z)) \cdot s'(z).
 \]
 Thus,
 \[
 |(r(s(z)))'| = |r'(s(z))| \cdot |s'(z)|.
 \]
 It follows that
 \begin{equation}{\label{4.5}}
  |(r(s(z)))'| \geq |r'(s(z))| \cdot \min_{z \in T_1} |s'(z)|.  \tag{26}
 \end{equation}
 From Lemma \ref{lem 2.2}, we obtain that 
 \[
 (r(s(z)))' \geq r'(s(z)) \cdot m \cdot \min_{z \in T_1} |s(z)|.
 \]
 Thus,
 \[
 (r(s(z)))' \geq mm' \cdot r'(s(z)). \tag{27}
 \]
 
 Therefore, it follows from (\ref{4.4}) that 
 \begin{equation*}
 	\left|{r}'\left( s(z)\right)  \right|\leq \frac{1}{2mm^{'}}\left[\left|{B}'(z) \right|-\frac{2\left[  \frac{mn}{2}- \left( \sum_{j=1}^{mt} \frac{1}{1+\left|b_{j} \right| }\right) \right]   \left(\left|r(s(z)) \right| -m^{*} \right)^{2}}{\left(\left \| r \circ s \right \|-m^{*} \right)^{2}}  \right] \\
 	 \times \left(\left \| r \circ s \right \|-m^{*} \right). 
 \end{equation*}
In the case where \( R(z) = 0 \), we find that \( (r(s(z)))' = 0 \).  \\
This implies that the above inequality is trivially satisfied.  
Therefore, inequality (18) holds for all \( z \in T_1 \). \\ 

Next, we demonstrate that equality is achieved for  
\[
r(s(z)) = \frac{(z + k)^{mt}}{(z - a)^{mn}},
\]  
where \( s(z) = z^m \) and  
\[
B(z) = \left(\frac{1 - az}{z - a}\right)^{mn}, \quad a > 1, \quad k \geq 1 \text{ at } z = 1.
\]  

First, observe that  
\[
\|r \circ s\| = \frac{(1 + k)^{mt}}{(a - 1)^{mn}} = |r(s(1))|,
\]  
and  
\[
m' = 1, \quad m^* = 0, \quad |B'(1)| = \frac{mn(a + 1)}{a - 1}.
\]  

For the derivative of \( r(s(z)) \), we have  
\[
r'(s(z)) = \frac{(z + k)^{mt}}{m z^{m - 1} (z - a)^{mn}} 
\left[\frac{mt}{z + k} + \frac{mn}{a - z}\right].
\]  

Evaluating at \( z = 1 \), we get  
\[
r'(s(1)) = \frac{(1 + k)^{mt}}{m(1)^{m - 1}(a - 1)^{mn}} 
\left(\frac{mt}{1 + k} + \frac{mn}{a - 1}\right) 
= \left(\frac{t}{1 + k} + \frac{n}{a - 1}\right)\|r \circ s\|.
\]  

The right-hand side inequality (\ref{thm2.1}) becomes  
\[
\frac{1}{2mm'}\left[|B'(1)| - \frac{(mn(1 + k) - 2mt)(|r(s(1))| - m^*)^2}{(1 + k)(\|r \circ s\| - m^*)^2}\right] (\|r \circ s\| - m^*).
\]  

Substituting the values, we find  
\[
\frac{1}{2m(1)}\left[\frac{mn(a + 1)}{a - 1} - \frac{mn(1 + k) - 2mt}{1 + k}\right] \times \|r \circ s\|.
\]  

Simplifying further,  
\[
\frac{1}{2}\left[\frac{2n}{a - 1} + \frac{2t}{1 + k}\right] \times \|r \circ s\| = |r'(s(1))|.
\]  

Thus, the bound is sharp and cannot be improved.

	\end{proof}
	\begin{proof}[\textbf{Proof of Theorem 2.2}]Suppose \( r(s(z)) \in \Re_{nm} \) has no zeros in \( D_{k}^{+} \), where \( k \leq 1 \).  
		Define \( m^{*} = \min_{\left| z \right| = k} \left| r(s(z)) \right| \). Then,  
		\[
		m^{*} \leq \left| r\left(s(z)\right) \right| \quad \text{for } z \in T_{k}.
		\]  
		
		If \( r\left(s(z)\right) \) has a zero on \( \left| z \right| = k \), then \( m^{*} = 0 \). Consequently, for any \( \alpha \) satisfying \( \left| \alpha \right| < 1 \),  
		\[
		r\left(s(z)\right) + \alpha m^{*} = r\left(r(z)\right).
		\]  
		
		On the other hand, if \( r\left(s(z)\right) \) has no zeros on \( \left| z \right| = k \), then for any \( \alpha \) with \( \left| \alpha \right| < 1 \),  
		\[
		\left| \alpha m^{*} \right| < \left| r\left(s(z)\right) \right| \quad \text{for } \left| z \right| = k.
		\]  
		
		By Rouche's theorem, it follows that \( R(z) = r\left(s(z)\right) + \alpha m^{*} \) and \( r\left(s(z)\right) \) have the same number of zeros in \( D_{k}^{-} \). Thus, for every \( \alpha \) with \( \left| \alpha \right| < 1 \),  
		\( R(z) \) has no zeros in \( D_{k}^{+} \).  
		
		If \( b_{1}, b_{2}, \dots, b_{mt} \) are the zeros of \( R(z) \) with \( mt \leq mn \),  
		and \( \left| b_{j} \right| \leq k \leq 1 \), we have  
		
	\begin{align*}
	\frac{z{R}'\left(z \right) }{R\left(z \right) }
	=  \frac{z({r(s(z))})'}{r(s(z)) }
	&=	 \frac{z({p(s(z))})' }{p(s(z)) } - \frac{z{w(s(z))}'}{w\left(s(z) \right) }\\
	&= \sum_{j=1}^{mt} \frac{z}{z - b_j} - \sum_{j=1}^{mn} \frac{z}{z - a_j}.
\end{align*}
Then
\begin{align*}
	Re\left(\frac{z{R}'\left(z \right) }{R\left(z \right) } \right) 
	= Re \left(\sum_{j=1}^{mt}\frac{z}{z-b_{j}} \right)  - Re\left(\sum_{j=1}^{mn} \frac{z}{z - a_j} \right) .
\end{align*}
\begin{align*}{\label{28}}
	Re\left(\frac{z{R}'\left(z \right) }{R\left(z \right) } \right) 
	&\geq \sum_{j=1}^{mt}\frac{1}{1+\left|b_{j} \right| }  + \sum_{j=1}^{mn} Re\left( \frac{1}{2}-\frac{z}{z - a_j}\right) -\frac{mn}{2}  \nonumber\\
	& = \sum_{j=1}^{mt}\frac{1}{1+\left|b_{j} \right| } -\frac{mn}{2}+ \sum_{j=1}^{mn} \frac{|a_j|^2 - 1}{2|z - a_j|^2}  \nonumber\\
	& = \sum_{j=1}^{mt}\frac{1}{1+\left|b_{j} \right| } -\frac{nm}{2}+\frac{|B^{'}(z)|}{2}. \tag{28} \nonumber\\
\end{align*}
Where $R\left(z \right)\neq 0, $ then
	\begin{align*}
		\left|\frac{{R}'\left(z \right) }{R\left( z\right) } \right|
		& = \left|\frac{z{R}'\left(z \right) }{R\left( z\right) } \right|\geq Re \left(\frac{z{R}'\left(z \right) }{R\left( z\right) }\right) 
		& \geq  \frac{|B^{'}(z)|}{2} -\frac{mn}{2} +\sum_{j=1}^{mt}\frac{1}{1+\left|b_{j} \right| }.	
	\end{align*}
	This implies that\\\\
	$ \left|{R}'\left(z \right)  \right|\geq \left[ \frac{|B^{'}(z)|}{2} -\frac{mn}{2} +\sum_{j=1}^{mt}\frac{1}{1+\left|b_{j} \right| }  \right]\left|R\left(z \right)  \right|, $ \text{for} $z \in T_{1}.$\\\\
	By replacing $R\left(z \right)=r\left( s(z)\right)+ \alpha m^{*} , $ we obtain that\\\\
	$ \left|{(r(s(z)))}'  \right|\geq \left[ \frac{|B^{'}(z)|}{2} -\frac{mn}{2} +\sum_{j=1}^{mt}\frac{1}{1+\left|b_{j} \right| }  \right]\left|r\left(s(z) \right) +\alpha m^{*} \right|$, \text{for} $z\in T_{1}.$\\\\
Observe that this inequality is trivially satisfied when \( R(z) = 0 \). Thus, it holds for all \( z \in T_{1} \). By appropriately selecting the argument of \( \alpha \) on the right-hand side of the inequality and noting that the left-hand side does not depend on \( \alpha \), we conclude that  

	$$ \left|({r(s(z))})'  \right|\geq \left[ \frac{|B^{'}(z)|}{2} -\frac{mn}{2} +\sum_{j=1}^{mt}\frac{1}{1+\left|b_{j} \right| }  \right]\left( \left|r\left(s(z) \right)\right|  +\left| \alpha\right|  m^{*}  \right),$$ for $z\in T_{1}.$\\
	Choose $\alpha$ such that $\left|\alpha \right|\rightarrow 1 $ for $z \in T_{1}$ that
	\begin{align*}
		\implies \left|{(r(s(z)))}'  \right|
		&\geq \left[ \frac{\left|{B}'\left(z \right)  \right| }{2} + \left(\sum_{j=1}^{mt}\frac{1}{1+\left|b_{j} \right| }-\frac{mn}{2} \right)  \right]\left( \left|r\left(s(z) \right)\right|  + m^{*} \right)\\
		&= \frac{1}{2}\left[ \left|{B}'\left(z \right)  \right| + 2\left(\sum_{j=1}^{mt}\frac{1}{1+\left|b_{j} \right| }-\frac{mn}{2} \right)  \right]\left( \left|r\left(s(z) \right)\right|  + m^{*} \right),
	\end{align*}
	which proves the Theorem \ref{thm2.2}.
\end{proof}

	\bibliography{ref}

\end{document}